\documentclass[a4paper,11pt]{amsart}
\usepackage{amssymb}
\author[A.~P.~Petravchuk]{Anatoliy P. Petravchuk}
\address{
Department of Algebra and Mathematical Logic\\ Faculty of Mechanics and Mathematics\\
Kyiv Taras Shevchenko University\\ Volodymyrska street, 01033
Kyiv, Ukraine} \email{aptr@univ.kiev.ua , apetrav@gmail.com}
\title[On pairs of commuting derivations of the polynomial ring in two variables ]
{On pairs of commuting derivations of the polynomial ring in two
variables}
\thanks{Partially supported by DFFD, Grant F28.1/026 }

\newcommand{\Ker}{\mathop{\mathrm{Ker}}}

\newtheorem{theorem}{Theorem}
\newtheorem{corollary}{Corollary}

\newtheorem{lemma}{Lemma}
\theoremstyle{definition}

\theoremstyle{remark}
\newtheorem{remark}{Remark}
\theoremstyle{problem}

\def\char{{\rm char\ }}
\def\div{{\rm div}}

\let\geq\geqslant
 
 \def\char{\rm char\,}
\begin{document}

\sloppy

\begin{abstract}

Let $k$ be an arbitrary field of characteristic zero, $k[x, y]$ be
the polynomial ring and $D$ a $k$-derivation of the ring $k[x,
y]$. Recall that a nonconstant polynomial $F\in k[x, y]$ is said
to be a Darboux polynomial of the derivation $D$ if $D(F)=\lambda
F$ for some polynomial $\lambda \in k[x, y]$. We prove that any
two linearly independent over the field $k$ commuting
$k$-derivations $D_{1}$ and $D_{2}$ of the ring $k[x, y]$ either
have a common Darboux polynomial, or   $D_{1}=D_{u_{1}}, \ \
D_{2}=D_{u_{2}}$ are Jacobian derivations i.e., $D_{i}(f)=\det
J(u_{i}, f)$ for every $f\in k[x, y], i=1, 2,$ where the
polynomials $u_{1}, u_{2}$ satisfy the condition $\det J(u_{1},
u_{2})=c\in k^{\star}.$ This statement about derivations  is an
analogue of the known fact from
 Linear Algebra about common eigenvectors of pairs of commuting
linear operators.

\end{abstract}

\maketitle

Let $k$ be an arbitrary field of characteristic zero, $k[x, y]$ be
the polynomial ring over $k.$
 One of the important questions in study of derivations of
polynomials rings  is the question about existence  Darboux
polynomials for a given $k$-derivation $D\in Der(k[x_{1}, \ldots ,
x_{n}])$ (recall that
 a polynomial $F\in k[x_{1}, \ldots , x_{n}]\setminus k$ is said to be
 a Darboux polynomial for a derivation $D$ if there
 exits $\lambda \in k[x_{1}, \ldots , x_{n}]$ such that
$D(F)=\lambda F$). There are many papers about  Darboux
polynomials of derivations (see, for example \cite{ Now1}, \cite
{Olla}).  It is known  that for every $n\geq 2$ the polynomial
ring $k[x_{1}, \ldots, x_{n}]$ has  a derivation without Darboux
polynomials (see, for example  \cite{Olla1}).

Since Darboux polynomials are in a sense eigenfunctions of the
linear operator $D$ on the vector space $k[x_{1}, \ldots , x_{n}]$
with polynomial eigenvalues, one can ask a natural question
whether given two commuting derivations $D_{1}, D_{2}$ of the
polynomial ring $k[x_{1}, \ldots , x_{n}]$ have a common Darboux
polynomial. If $\char k>0,$ then the answer to this question is
positive, because every nonconstant polynomial of the form
$f(x_{1}^{p}, \ldots , x_{n}^{p})$ is a Darboux polynomial (with
zero eigenvalue) for any derivation $D\in Der(k[x_{1}, \ldots ,
x_{n}]).$ But in zero characteristic the answer is negative: for
example, the derivations $D_{1}=\frac{\partial}{\partial x}$ and
$D_{2}=\frac{\partial}{\partial y}$ of the polynomial ring $k[x,
y]$ commute, but they have not any common Darboux polynomials.
Besides, the case is possible when two derivations $D_{1}$ and
$D_{2}$ are linearly dependent over the field $k$ (and
consequently commuting) and they have not any Darboux polynomials.

The main result of this paper shows that  with the exception of
above mentioned examples any two commuting derivations of the
polynomial ring $k[x, y]$ always have  a common Darboux
polynomial.

The notations in the paper are standard (see, for example,
\cite{Now1}). For an arbitrary polynomial $u\in k[x, y]$ we denote
by $D_{u}$ the Jacobian derivation determined by the rule
$D_{u}(f)=\det J(u, f)$ for $f\in k[x, y]$, where $J(u, f)$ is the
Jacoby matrix of the polynomials $u, f.$ The Lie algebra of all
$k$-derivations of the polynomial ring $k[x, y]$ will be denoted
by $W_{2}(k)$. Since every derivation $D\in W_{2}(k)$ can be
uniquely written  in the form $D=P(x, y)\frac{\partial}{\partial
x}+ Q(x, y)\frac{\partial}{\partial y}$ where $P(x, y)=D(x), Q(x,
y)=D(y)$ we see that $W_{2}(k)$ is a free $k[x, y]$-module of rank
$2.$ For a derivation $D=P(x, y)\frac{\partial}{\partial x}+ Q(x,
y)\frac{\partial}{\partial y}$ we denote by  $\div D$ the
divergence $\frac{\partial P}{\partial x}+\frac{\partial
Q}{\partial y}$  of $D$.

 A derivation $D=P(x, y)\frac{\partial}{\partial x}+
Q(x, y)\frac{\partial}{\partial y}$ will be called reduced if $P,
Q$ are coprime as polynomials from $k[x, y].$ Every derivation
$D\in W_{2}(k)$ can be uniquely extended to a derivation of the
field $k(x, y)$ of all rational functions in two variables if we
put $D(\frac{f}{g})=\frac{D(f)g-fD(g)}{g^{2}}$.

\begin{lemma}\label{dependence}

Let $D_{1}, D_{2}\in W_{2}(k)$ and $D_{1}, D_{2}$ be linearly
dependent over the ring $k[x, y]$. Then there exist a reduced
derivation $D_{0}$ and polynomials $f, g\in k[x, y]$ such that
$D_{1}=fD_{0}, D_{2}=gD_{0}$. If $[D_{1}, D_{2}]=0,$ then
$D_{0}(\frac{f}{g})=0.$
\end{lemma}
\begin{proof}
Let $$D_{1}=P_{1}\frac{\partial}{\partial
x}+Q_{1}\frac{\partial}{\partial y}, \
D_{2}=P_{2}\frac{\partial}{\partial
x}+Q_{2}\frac{\partial}{\partial y}$$
 and $\mu _{1}=\gcd (P_{1},
Q_{1}), \ \mu _{2}=\gcd (P_{2}, Q_{2}).$ Denote by ${\overline
D}_{1}=D_{1}/\mu _{1}, {\overline D}_{2}=D_{2}/\mu _{2}$ the
corresponding reduced derivations. Since $\alpha D_{1}+\beta
D_{2}=0$ by conditions of Lemma for some $\alpha , \beta\in k[x,
y]$, we also have  $\alpha \mu _{1}{\overline D}_{1}+\beta \mu
_{2}{\overline D}_{2}=0.$ Denote $\lambda _{1}=\alpha \mu _{1},
\lambda _{2}=\beta \mu _{2}.$  It follows from the equality
$\lambda _{1}{\overline D}_{1}+\lambda _{2}{\overline D}_{2}=0$
that
$$\lambda _{1}{\overline P}_{1}+\lambda _{2}{\overline P}_{2}=0, \
\ \lambda _{1}{\overline Q}_{1}+\lambda _{2}{\overline Q}_{2}=0,
\eqno (1)$$ where ${\overline P}_{i}=P_{i}/\mu _{i}, \ \
{\overline Q}_{i}=P_{i}/\mu _{i}, i=1, 2.$ From the equalities (1)
it follows that
$$\frac{\lambda _{1}}{\lambda _{2}}=- \frac{{\overline P}_{1}}{{\overline
P}_{2}}=- \frac{{\overline Q}_{1}}{{\overline Q}_{2}}$$ and
therefore ${\overline P}_{1}{\overline Q}_{2}={\overline
P}_{2}{\overline Q}_{1}.$ As $\gcd ({\overline P}_{1}, {\overline
Q}_{1})=1,$ we have ${\overline P}_{1}=\lambda {\overline P}_{2}$
and ${\overline Q}_{1}=\lambda {\overline Q}_{2}$ for some
$\lambda \in k[x, y].$ One can analogously show that ${\overline
P}_{2}=\nu {\overline P}_{1}$ and ${\overline Q}_{2}=\nu
{\overline Q}_{1}$ for some $\nu \in k[x, y].$ Then it is easily
shown that ${\overline D}_{2}=c{\overline D}_{1}$ for some $c\in
k^{\star}.$ Denote $D_{0}={\overline D}_{1}.$ Then $D_{1}=fD_{0},
\ D_{2}=gD_{0}$ for some polynomials $f, g\in k[x, y].$

Let now $[D_{1}, D_{2}]=0.$ It can  easily be checked that
$$[D_{1}, D_{2}]=[fD_{0}, gD_{0}]=(D_{0}(f)g-fD_{0}(g))D_{0}=0.$$
This yields the equality $D_{0}(f)g-fD_{0}(g)=0.$
 Consequently,
 $$D_{0}(\frac{f}{g})=\frac{D_{0}(f)g-fD_{0}(g)}{g^{2}}=0.$$
\end{proof}

The next result   is  a very special case of  Theorem 2.4 from
\cite{Buchstaber}, where polynomial Lie algebras over the field
$\mathbb C$ were studied. We give here an elementary proof of this
special result  over an arbitrary field.

\begin{lemma}\label{Buchstaber}
Let $D_{1}=P_{1}\frac{\partial}{\partial
x}+Q_{1}\frac{\partial}{\partial y},
D_{2}=P_{2}\frac{\partial}{\partial
x}+Q_{2}\frac{\partial}{\partial y}$ be commuting derivations of
the polynomial ring $k[x, y]$ over an arbitrary field $k$. Then
 for the determinant $\Delta =\left|
\begin{array}{cc}
  P_{1} & Q_{1} \\
  P_{2} & Q_{2} \\
\end{array}%
\right |$ it holds $ D_{i}(\Delta )=\Delta \cdot \div{D_{i}}, \
i=1, 2.$
\end{lemma}
\begin{proof}
Prove, for example, that $D_{1}(\Delta )=\Delta \cdot
\div{D_{1}}=\Delta (P_{1x}^{'}+Q_{1y}^{'}).$ We have equalities
$$D_{1}(\Delta )=D_{1}\left ( \left|
\begin{array}{cc}
  P_{1} & Q_{1} \\
  P_{2} & Q_{2} \\
\end{array}%
\right | \right )= \left |
\begin{array}{cc}
  D_{1}(P_{1}) & D_{1}(Q_{1}) \\
  P_{2} & Q_{2} \\
\end{array}%
\right |+ \left|
\begin{array}{cc}
  P_{1} & Q_{1} \\
  D_{1}(P_{2}) & D_{1}(Q_{2})\\
\end{array}%
\right | .$$ Since $[D_{1}, D_{2}]=0$, we have
$D_{1}(P_{2})=D_{2}(P_{1})$ and $D_{1}(Q_{2})=D_{2}(Q_{1})$ ( this
can be shown by straightforward check). Therefore
$$D_{1}(\Delta )=\left |
\begin{array}{cc}
  D_{1}(P_{1}) & D_{1}(Q_{1}) \\
  P_{2} & Q_{2} \\
\end{array}%
\right |+ \left|
\begin{array}{cc}
  P_{1} & Q_{1} \\
  D_{2}(P_{1}) & D_{2}(Q_{1})\\
\end{array}%
\right | =$$
$$=Q_{2}D_{1}(P_{1})-P_{2}D_{1}(Q_{1})+P_{1}D_{2}(Q_{1})-Q_{1}D_{2}(P_{1}).$$
After substituting $D_{i}(P_{1})$ and $D_{i}(Q_{1})$ in the last
equality we obtain
$$D_{1}(\Delta
)=Q_{2}(P_{1}P_{1x}^{'}+Q_{1}P_{1y}^{'})-P_{2}(Q_{1}P_{1x}^{'}+Q_{1}Q_{1y}^{'})+
P_{1}(P_{2}Q_{1x}^{'}+Q_{2}Q_{1y}^{'})-$$
$$-Q_{1}(P_{2}P_{1x}^{'}+Q_{2}P_{1y}^{'}).$$
 After opening the
brackets and cancellation we get
$$D_{1}(\Delta
)=P_{1}Q_{2}P_{1x}^{'}-Q_{1}P_{2}P_{1x}^{'}+P_{1}P_{2}Q_{1y}^{'}-Q_{1}P_{2}Q_{1y}^{'}=$$
$$=\left |
\begin{array}{cc}
  P_{1} & Q_{1} \\
  P_{2} & Q_{2} \\
\end{array} \right | P_{1x}^{'}+\left|
\begin{array}{cc}
  P_{1} & Q_{1} \\
  P_{2} & Q_{2} \\
\end{array} \right | Q_{1y}^{'}=\Delta \cdot  \div{D_{1}}.$$
One can analogously show that $D_{2}(\Delta )=\Delta \cdot
\div{D_{2}}.$

\end{proof}

\begin{lemma}\label{Jacob}
Let $D_{u_{1}}, D_{u_{2}}$ be Jacobian derivations of the ring
$k[x, y]$ such that $\det J(u_{1}, u_{2})=c\in k^{\star}.$ Then
the derivations $D_{u_{1}}, D_{u_{2}}$ do not have any common
Darboux polynomials.
\end{lemma}

\begin{proof}
Let the statement of Lemma be false and $F$ be a common Darboux
polynomial for $D_{u_{1}}$ and $ D_{u_{2}}$. By Proposition 2.2.1
from \cite{Now1}
 every nonconstant divisor of $F$ is a Darboux polynomial for
the both derivation $D_{1}$ and $D_{2}$. Therefore without loss of
generality one can assume that the polynomial $F$ is irreducible.
By Lemma 2.1 from \cite{St} there exist elements $c_{1}, \
c_{2}\in k$ such that $u_{1}-c_{1}$ and $u_{2}-c_{2}$ are
divisible by $F.$ Let $u_{1}-c_{1}=\lambda _{1}F$ and
$u_{2}-c_{2}=\lambda _{2}F$ for some $\lambda _{1}, \lambda
_{2}\in k[x, y].$ Then we have
$$\det J(u_{1}, u_{2})=\det J(u_{1}-c_{1}, u_{2}-c_{2})=\det
J(\lambda _{1}F, \lambda _{2}F).$$ But it is easily shown that
$$\det J(\lambda _{1}F, \lambda _{2}F)=F^{2}\det J(\lambda _{1},
\lambda _{2})+\lambda _{1}F\det J(F, \lambda _{2})+\lambda
_{2}F\det J(\lambda _{1}, F)$$ and consequently the polynomial
$\det J(u_{1}, u_{2})$ is divisible by $F.$ This is impossible
because $\det J(u_{1}, u_{2})=c\in k^{\star}.$ This contradiction
concludes the proof.
\end{proof}

\

\begin{theorem}\label{main}
Let $k[x, y]$ be the polynomial ring over an arbitrary field $k$
of characteristic zero and $D_{1}, D_{2}$ be linearly independent
over the field $k$ derivations of $k[x, y]$. If $[D_{1}, D_{2}]=0$
then either $D_{1}$ and $D_{2}$ have a common Darboux polynomial,
or $D_{1}=D_{u_{1}}, D_{2}=D_{u_{2}}$ are Jacobian derivations,
where $u_{1}, u_{2}$ are polynomials with $\det J(u_{1},
u_{2})=c\in k^{\star}.$

\end{theorem}

\begin{proof}
Let $D_{1}=P_{1}\frac{\partial}{\partial
x}+Q_{1}\frac{\partial}{\partial y}$ and $
D_{2}=P_{2}\frac{\partial}{\partial
x}+Q_{2}\frac{\partial}{\partial y}$. Denote  $\Delta =\left|
\begin{array}{cc}
  P_{1} & Q_{1} \\
  P_{2} & Q_{2} \\
\end{array}%
\right |$ and consider at first the case when $\Delta =0.$ Then
the rows of the matrix $\left (
\begin{array}{cc}
  P_{1} & Q_{1} \\
  P_{2} & Q_{2} \\
\end{array}%
\right )$ are linearly dependent over the ring $k[x, y]$ and hence
$\alpha D_{1}+\beta D_{2}=0$ for some $\alpha , \beta \in k[x,
y].$ By Lemma~\ref{dependence} there exists a reduced derivation
$D_{0}$ such that $D_{1}=fD_{0}, \ D_{2}=gD_{0}$ for some $f, g\in
k.$ Since $[D_{1}, D_{2}]=0$, it follows from the same Lemma that
$D_{0}(\frac{f}{g})=0.$  Note that the rational function
$\frac{f}{g}$ is nonconstant because otherwise the derivations
$D_{1}$ and $D_{2}$ were linearly dependent over $k$ which is
impossible.

Consider $D_{0}$ as a derivation of the field of rational
functions $k(x, y).$ Then $\frac{f}{g}\in \Ker D_{0}$ and $\Ker
(D_{0})$  is an algebraically closed subfield of the field $k(x,
y)$ (see for example \cite{Now1}, Lemma 2.1).

 By Corollary 2 from \cite{Petien2} either $\Ker D_{0}=k(p)$ for some
irreducible polynomial $p$ or $\Ker D_{0}=k(\frac{p}{q})$ where
$p, q$ are algebraically independent irreducible polynomials. In
the first case the polynomial $p$ is a common Darboux polynomial
for $D_{1}$ and $D_{2}$ with eigenvalue $\lambda =0.$

In the second case $D_{0}(\frac{p}{q})=0$ and by Proposition 2.2.2
from \cite{Now1} $D_{0}(p)=\mu p, \ D_{0}(q)=\mu q$ for some $\mu
\in k[x, y].$ It is obvious that $p$ is a common Darboux
polynomial for $D_{1}$ and $D_{2}$ (as well as  the polynomial
$q$).

Consider now the case $\Delta =\left|
\begin{array}{cc}
  P_{1} & Q_{1} \\
  P_{2} & Q_{2} \\
\end{array}%
\right |\not =0.$ If $\Delta \not= {\rm const}$, then by
Lemma~\ref{Buchstaber} $\Delta $ is a common Darboux polynomial
for $D_{1}$ and $D_{2}.$ Let now $\Delta =c\in k^{\star}.$ Then by
Proposition 2.2.4 from \cite{Essen} the derivations $D_{1},$
$D_{2}$ form a basis of the $k[x, y]$-module $W_{2}(k)$ and are
both the Jacobian derivations $D_{1}=D_{u_{1}}$ and
$D_{2}=D_{u_{2}}$ with $\det J(u_{1}, u_{2})=\Delta \in
k^{\star}.$ This completes the proof.
\end{proof}
\begin{remark}
By Lemma~\ref{Jacob} any two commuting derivations
$D_{1}=D_{u_{1}}$ and $D_{2}=D_{u_{2}}$ with $\det J(u_{1},
u_{2})=c\in k^{\star}$ do not have any common Darboux polynomials.
\end{remark}

\begin{corollary}
Let $D$ be a derivation of the ring $k[x, y]$ with a Darboux
polynomial. If $D_{1}\in W_{2}(k)$ and $[D, D_{1}]=0$, then the
derivation $D_{1}$ also has a Darboux polynomial.
\end{corollary}
\begin{proof}
If $D$ and $D_{1}$ do not have any common Darboux polynomial, then
 by Theorem~\ref{main} either $D$ and $D_{1}$ are linearly
 dependent over $k$ or $D=D_{u}$ and $D_{1}=D_{u_{1}}$ for some
 polynomials $u, u_{1}\in k[x, y]$ with $\det J(u, u_{1})=c\in
 k^{\star}.$ The first case is impossible  because  the derivation $D_{1}$ has the same
 Darboux polynomial as $D$, in the second case $D_{1}=D_{u_{1}}$
 has the polynomial Darboux $u_{1}$ with the eigenvalue $\lambda
 =0.$
\end{proof}


\begin{thebibliography}{99}


\bibitem{Bodin}
 Bodin A.
 {\it Reducibility of rational functions in several variables},
 Israel J. Math. 164 (2008), 333--347.

\bibitem{Buchstaber}
Buchstaber V.M., Leykin D.V.
 {\it Polynomial Lie algebras},
 Functional Analysis and its Applications, 36, no.4
(2002),   267-280.



\bibitem{Now1}
Nowicki A., {\it Polynomial derivations and their rings of
constants}, N. Copernicus University Press: Toru${\acute {\rm
n}}$, 1994.

\bibitem{Now2}
Nowicki A. and Nagata M. {\it  Rings of constants for
{$k$}-derivations in $k[x\sb 1,\dots,x\sb
  n]$} J. Math. Kyoto Univ.,  28(1) (1988), 111--118.

\bibitem{Olla1}
Ollagnier J.M. and Nowicki A. {\it Derivations of polynomial
algebras without Darboux polynomials}, Journ. of Pure and Appl.
Algebra, 212 (2008), 1626-1631.

\bibitem{Olla}
 Ollagnier J.~M. {\it  Algebraic closure of a rational function},
 Qualitative theory of dynamical systems, 5(2) 2004, 285--300.


\bibitem{Petien2}
 Petravchuk A.P. and Iena O.G.  {\it  On closed rational
functions in several variables},  Algebra Discrete Math., (2)
2007, 115--124.

\bibitem{St}
Stein Y. \textit{The total reducibility order of a polynomial in
two variables}, Israel J. Math. {\bf 68} (1989), 109--122.

\bibitem{Essen}
van den Essen, A., Polynomial automorphisms and the Jacobian
Conjecture, Progress in Mathematics, 190, Birkhauser Verlag, 2000.

\end{thebibliography}
\end{document}